\documentclass{amsart}

\usepackage{amsmath, amsfonts, mathrsfs}
\usepackage{amssymb,latexsym, fontenc}
\usepackage{enumerate}
\usepackage{hyperref}
\usepackage{bbm}

\usepackage[top=3.1cm, bottom=3.1cm, left=3.1cm, right=3.1cm, includefoot, heightrounded]{geometry}

\newtheorem{theorem}{Theorem}[section]

\newtheorem{lemma}[theorem]{Lemma}

\theoremstyle{definition}
\newtheorem{definition}[theorem]{Definition}

\numberwithin{equation}{section}

\begin{document}

\title[A full characterization of Pettis Integrable Multifunctions]
{A full characterization of Pettis Integrable Multifunctions in Fr\'{e}chet Spaces}

\author{ Sokol  Kaliaj }

\address{
Mathematics Department, 
Science Natural Faculty, 
"Aleksand\"{e}r Xhuvani" University,
Elbasan, 
Albania.
}

\email{printimesokol@gmail.com}

\thanks{}



\begin{abstract}
In this paper, 
we present a full characterization of Pettis integrability in $cwk(X)$ of a scalarly integrable multifunction 
$\Gamma : \Omega \to cwk(X)$ 
in terms of weakly compact generated (\text{WCG}) subspaces of a Fr\'{e}chet space $X$, 
where $(\Omega, \Sigma, \mu)$ is a complete finite measure space 
and 
$cwk(X)$ is the family of all nonempty convex weakly compact subsets of the Fr\'{e}chet space $X$.  
A convergence theorem for Pettis integrable multifunctions is proved. 
This result yields  a convergence theorem for Pettis integrable fuzzy maps defined on $\Omega$ and taking values 
in the generalized fuzzy number hyperspace 
$\mathcal{F}_{cwk}(X)$. 
\end{abstract}

\maketitle

\section{Introduction and Preliminaries}

There are characterizations of Pettis integrability in $bcc(X)$ and $cwk(X)$ 
of multifunctions defined on $\Omega$ and taking values in the hyperspace $cc(X)$ 
(see \cite[Lemma 4.3, Lemma 4.5, Theorem 4.6]{KAL_Pettis_LCS}), 
where $cc(X) (bcc(X), cwk(X))$ is the family of all nonempty 
convex closed  (bounded  convex closed, convex weakly compact) subsets of a 
Hausdorff locally convex topological vector space $X$.  
It is known a full characterization 
of Pettis integrability in $cwk(X)$ 
of a scalarly integrable multifunction 
$\Gamma : \Omega \to cwk(X)$ 
in terms of \text{WCG} subspaces of a Banach space $X$, 
see \cite[Theorem 2.5]{MUS_Pettis_MultiF}. 
In this paper, 
we first present  
a full characterization of Pettis integrability in $cwk(X)$ of a scalarly integrable multifunction 
$\Gamma : \Omega \to cwk(X)$ 
in terms of weakly compact generated (\text{WCG}) subspaces of a Fr\'{e}chet space $X$, 
 Theorem \ref{T_Pettis_WCG}.
Then, a convergence theorem 
for Pettis integrable multifunctions is proved, Theorem \ref{CT_Pettis_MultiF}.

It is known the close relation between 
Pettis integrable fuzzy maps and Pettis integrable multifunctions. 
Thus, for example,   
the relation between 
Radon-Nikodym derivative in terms of Pettis integral 
of fuzzy number multimeasures 
and of multimeasures 
(see \cite{BAN, STOJ_RNT, WU_RNT, ZHOU, KAL_GFNM, KAL_FNM, KAL_FNM_Frechet}), 
or the relation between full characterizations 
of Pettis integrable fuzzy maps (\cite[Theorem 5.5]{KAL_Pettis_LCS}) 
and of Pettis integrable multifunctions (\cite[Theorem 4.6]{KAL_Pettis_LCS}). 
Here, a convergence theorem for Pettis integrable fuzzy maps (Theorem \ref{CT_Pettis_FuzzyMaps}) is proved by 
the convergence theorem for Pettis integrable multifunctions (Theorem \ref{CT_Pettis_MultiF}).

Throughout this paper, $(X,\tau)$ is a Fr\'{e}chet spaces and $(p_{i})_{i \in \mathbb{N}}$ 
is an increasing sequence of seminorms defining the topology $\tau$.
For any $i \in \mathbb{N}$,  
$X_{i}$ is the quotient vector space $X/p_{i}^{-1}(0)$, 
$\varphi_{i}: X \to X_{i}$ is 
the canonical quotient map,  
$(X_{i},||\cdot||_{i})$ is the quotient normed space 
($||\varphi_{i}(x))||_{i} = p_{i}(x)$ for all $x \in X$), 
and $\overline{X}_{i}$ is the completion of $(X_{i},||\cdot||_{i})$.
By $X^{*}$ is denoted the topological dual of $X$, 
by $X^{*}_{i}$ is denoted the topological dual of the quotient normed space  
$(X_{i},||\cdot||_{i})$, 
by $B_{X_{i}^{*}}$ is denoted the closed unit ball in $X^{*}_{i}$, 
and by $X^{*}_{p_{i}}$ is denoted the space of all linear functionals which are continuous with respect to $p_{i}$, 
i.e., 
\begin{equation*}
\begin{split}
X^{*}_{i} = \mathcal{L}((X_{i},||\cdot||_{i}), \mathbb{R})
\quad
\text{and}
\quad
X^{*}_{p_{i}} = \mathcal{L}((X,p_{i}), \mathbb{R})
\qquad
(i \in \mathbb{N}),
\end{split}
\end{equation*}
where 
$\mathcal{L}((\mathfrak{X}, \mathfrak{T}), \mathbb{R})$ 
is the family of all continuous linear functional defined on 
a locally convex space $(\mathfrak{X}, \mathfrak{T})$ 
and taking values in $\mathbb{R}$. 
It is clear that $X^{*}_{p_{i}} \subset X^{*}_{p_{i+1}}$ for all $i \in \mathbb{N}$, and 
$X^{*} = \bigcup_{i=1}^{+\infty} X^{*}_{p_{i}}$.

By $\mathcal{P}(X)$ is denoted the family of all non-empty subsets of $X$. 
For each $C \in \mathcal{P}(X)$ 
the \textit{support} of $C$ is defined by equality
\begin{equation*}
\begin{split}
\sigma(x^{*}, C) = 
\sup \{ x^{*}(x) : x \in C \}
\qquad (x^{*} \in X^{*}).
\end{split}
\end{equation*}

The \textit{Hausdorff semimetric} associated with the seminorm $p_{i}$ is defined as follows 
\begin{equation*}
\begin{split}
H_{p_{i}}(A,B) = \max\{ e_{p_{i}}(A,B), e_{p_{i}}(B,A) \},
\quad
(A, B \in \mathcal{P}(X)), 
\end{split}
\end{equation*}
where 
\begin{equation*}
\begin{split}
e_{p_{i}}(A,B) = \sup_{a \in A} p_{i}(a,B), 
\quad  
p_{i}(a,B) = \inf \{p_{i}(a-b) : b \in B \}.
\end{split}
\end{equation*}
Note that if $h_{i}$ is defined as follows
\begin{equation*}
\begin{split}
h_{i}(A_{i},B_{i}) = 
\max \left \{ e_{i}(A_{i},B_{i}), e_{i}(B_{i},A_{i}) \right \},
\end{split}
\end{equation*}
where 
$A_{i} = \varphi_{i}(A)$, $B_{i} = \varphi_{i}(B)$, $(A, B \in \mathcal{P}(X))$, 
and 
\begin{equation*}
\begin{split}
e_{i}(A_{i},B_{i}) = \sup_{a \in A}  \inf \{||\varphi_{i}(a)-\varphi_{i}(b)||_{i} : b \in B \},
\end{split}
\end{equation*}
then
\begin{equation}\label{eq_Hausdorff_Xi}
\begin{split}
H_{p_{i}}(A,B) = h_{i}(A_{i},B_{i})
\qquad
(A, B \in \mathcal{P}(X), i \in \mathbb{N}).
\end{split}
\end{equation}
Given a seminorm $p_{i}$ we write 
$U_{p_{i}} = \{x \in X : p_{i}(x) \leq 1\}$ and denote by $U_{p_{i}}^{o}$ the absolute polar of $U_{p_{i}}$, 
i.e., $U_{p_{i}}^{o} = \{x^{*} \in X^{*} : \sup_{x \in U_{p_{i}}} |x^{*}(x)| \leq 1\}$.


If  $(B_{n})$ is a sequence with terms in $\mathcal{P}(X)$, 
the symbol $\sum_{n=1}^{+\infty} B_{n}$ 
denotes a formal series. 
The series $\sum_{n=1}^{+\infty} B_{n}$ is said to be 
\textit{unconditionally convergent} 
provided that for every choice $b_{n} \in B_{n}, n \in \mathbb{N}$, 
the series $\sum_{n=1}^{+\infty} b_{n}$ is unconditionally convergent in $X$. 
In this case we define
\begin{equation}\label{eqD_Multimeasures}
\sum_{n=1}^{+\infty} B_{n} = 
\left \{
\sum_{n=1}^{+\infty} b_{n}: b_{n} \in B_{n} \text{ for all }n \in \mathbb{N}
\right \}.
\end{equation}

\begin{definition}\label{D_Multimeasures}
A map $M: \Sigma \to bcc(X)$ is called a \textit{multimeasure},  
if the map $\sigma(x^{*}, M(\cdot))$ is a countable additive real valued measure for every 
$x^{*} \in X^{*}$. 
The map $M$ is said to be a \textit{strong multimeasure}, if
\begin{itemize}
\item[(i)]
$M(\emptyset) = \{\theta\}$, 
\item[(ii)]
for each sequence $(A_{n})$ of pairwise disjoint members of $\Sigma$, 
we have $M(\bigcup_{n} A_{n}) = \sum_{n} M(A_{n})$, 
where  $\sum_{n} M(A_{n})$ is defined by \eqref{eqD_Multimeasures}.  
\end{itemize}
It is clear that each strong multimeasure is a multimeasure. 
If $M$ is a multimeasure taking values in the hyperspace $cwk(X)$, then the converse is also true, 
see \cite[Theorem 2.10]{KAL_REND}. 
We say that $M$ is an \textit{$h$-multimeasure} 
if $M$ is countable additive with respect to the Hausdorff semimetric $H_{p_{i}}$ 
for every $i \in \mathbb{N}$.  
Thanks to \eqref{eq_Hausdorff_Xi}, 
$M$ is an $h$-multimeasure if and only if 
$M_{i}$ is an $h_{i}$-multimeasure for every $i \in \mathbb{N}$, 
where 
$M_{i}(A) = \varphi_{i}(M(A)) = \{\varphi_{i}(x): x \in M(A)\}$, 
($A \in \Sigma, i \in \mathbb{N}$).
We will show that if  $M: \Sigma \to cwk(X)$ is a multimeasure, 
then $M$ is an $h$-multimeasure, see Lemma \ref{L_Multi_HausdorffM}. 
\end{definition}


\begin{definition}\label{D_Pettis_Multi}
A multifunction $\Gamma: \Omega \to cc(X)$ is said to be 
\textit{scalarly measurable} (\textit{scalarly integrable}) 
if the map $\sigma(x^{*}, \Gamma(\cdot))$ is measurable (Lebesgue integrable) 
whenever $x^{*} \in X^{*}$. 
The multifunction $\Gamma$ is said to be
\textit{Pettis integrable} in $bcc(X)$ ($cwk(X)$),  
if 
\begin{itemize}
\item
$\Gamma$ is scalarly integrable,
\item
for each $A \in \Sigma$ there exists $M(A) \in bcc(X)$ ($M(A) \in cwk(X)$) such that
\begin{equation*}
\begin{split}
\sigma
\left ( 
x^{*}, M(A) 
\right )
=
\int_{A} \sigma(x^{*}, \Gamma(t)) d \mu
\quad
\text{for every }x^{*} \in X^{*}.
\end{split}
\end{equation*} 
\end{itemize}
We call $M(A)$ the \textit{Pettis integral} of $\Gamma$ over $A$ and 
set $(P)\int_{A} \Gamma(t) d \mu = M(A)$. 
It is clear that the map 
$M$ defined by $M(A) = (P)\int_{A} \Gamma(t) d \mu$ for all $A \in \Sigma$, 
is a multimeasure. 
\end{definition} 
The Pettis integral for multifunctions was first considered by Castaing and Valadier 
\cite[Chapter V]{CAST}. 
For the general theory of Pettis integrable multifunctions  
we use the papers 
of H.~Ziat, K.~El Amri and C.~Hess \cite{EL, ZIAT}, 
B.~Cascales, V.~Kadets, and J.~Rodr\'{i}guez \cite{CASC_SelectorsA, CASC_SelectorsB} and 
K.~Musial \cite{MUS_Pettis_MultiF, MUS_Approximation, MUS_Conjugate_Banach_Spaces,  MUS_RNT_Rich_Multimeasures_I, MUS_RNT_Rich_Multimeasures_II}. 
A presentation of the general theory of Pettis integration of vector functions can be
found in \cite{TAL} and \cite{MUS_Pettis_Integral_Functions}.

\begin{definition}\label{D_Uniformly_Integrable}
A family $\mathcal{F} \subset L_{1}(\mu)$ is said to be 
\textit{uniformly integrable} if $\mathcal{F}$ is bounded in $(L_{1}(\mu), ||\cdot||_{1})$ 
and for each $\varepsilon >0$ there exists $\delta >0$ such that 
for each $A \in \Sigma$, 
\begin{equation*}
\begin{split}
\mu(A) < \delta 
\Rightarrow
\sup_{f \in \mathcal{F}} \int_{A} |f(t)| d \mu < \varepsilon.
\end{split}
\end{equation*}
\end{definition}


\begin{definition}\label{D_Sublinear_Map}
A function $\varphi : X^{*} \to L_{1}(\mu)$ is said to be 
\begin{itemize}
\item
\textit{positively homogeneous} 
if $\varphi(\alpha x^{*}) = \alpha \varphi(x)$ for every 
$\alpha \geq 0$ and $x^{*} \in X^{*}$,
\item
\textit{sub-additive} 
if $\varphi(x^{*} +y^{*}) \leq \varphi(x^{*}) + \varphi(y^{*})$ for all pairs 
$(x^{*}, y^{*}) \in X^{*} \times X^{*}$, where
\begin{equation*}
\begin{split}
f \leq g \Leftrightarrow f(t) \leq g(t)\text{ at almost all }t \in \Omega
\quad
(f,g \in L_{1}(\mu)),
\end{split}
\end{equation*}
\item
\textit{sublinear} 
if it is both positively homogeneous and sub-additive. 
\end{itemize}
Given a multifunction $\Gamma : \Omega \to cc(X)$ scalarly integrable 
we can define a sublinear map 
\begin{equation*}
\begin{split}
\varphi_{\Gamma} : X^{*} \to L_{1}(\mu), 
\quad
\varphi_{\Gamma}(x^{*}) = \sigma(x^{*}, \Gamma(\cdot)). 
\end{split}
\end{equation*}
We say that $\varphi_{\Gamma}$ is \textit{weakly compact} if 
$\varphi_{\Gamma}(U_{p_{i}}^{o})$ is a relatively weakly compact subset of $L_{1}(\mu)$ 
for each $i \in \mathbb{N}$. 
\end{definition}


Thanks to James's theorem (\cite[Theorem 6]{JAMES}), 
the following one is proved in the same way as \cite[Proposition 2.4]{KAN}.

\begin{theorem}\label{T_CWK}
Let $M : \Sigma \to cwk(X)$ be a multimeasure. 
Then $M(\Sigma) = \bigcup_{E \in \Sigma} M(E)$ is a relatively weakly compact subset of $X$.
\end{theorem}

\begin{definition}\label{D_DeterminedMultiF_By_WCG}
A vector subspace $Y \subset X$ 
is said to be \textit{Weakly Compact Generated (WCG) subspace} of $X$ if there exists 
a weakly compact set 
$K \subset X$ such that $Y$ is 
the closed linear space generated by $K$,  
i.e., $Y = \overline{\text{sp}(K)}$. 
A multifunction $\Gamma: \Omega \to cc(X)$ is said to be \textit{determined} by 
a subspace $Y \subset X$
if 
\begin{equation*} 
\begin{split}
x^{*} \in Y^{\perp}
\Rightarrow 
\sigma(x^{*}, \Gamma(t)) = 0
~\text{ at almost all }t \in \Omega,
\end{split}
\end{equation*} 
where the exceptional sets depend on $x^{*}$, and 
$Y^{\perp} = \{x^{*} \in X^{*}: (\forall y \in Y)[x^{*}(y) = 0]\}$.  
\end{definition}

\begin{definition}\label{D_Fuzzy_Number}
A map $u: X \to [0,1]$ is said to be a \textit{fuzzy set} on $X$
(see Zadeh \cite{ZAD}).
For each $\alpha \in (0,1]$, denote by 
\begin{equation*}
\begin{split}
[u]^{\alpha} = \{x \in X: u(x) \geq \alpha\}
\end{split}
\end{equation*}
the $\alpha$-level set of $u$. 
We say that the fuzzy set $u$ on $X$ is a \textit{generalized fuzzy number} on $X$, 
if $[u]^{\alpha} \in bcc(X)$ for all $\alpha \in (0,1]$.  
We denote by $\mathcal{F}_{bcc}(X)$ the set of all generalized fuzzy numbers 
on $X$. 
In particular, 
by $\mathcal{F}_{cwk}(X)$ is denoted the set of all generalized fuzzy numbers 
$u \in \mathcal{F}_{bcc}(X)$ such that 
$[u]^{\alpha} \in cwk(X)$  for all $\alpha \in (0,1]$. 
\end{definition}

\begin{definition}\label{D_Fuzzy_Pettis}
Given a fuzzy map $\Gamma : \Omega \to \mathcal{F}_{cwk}(X)$ 
we can define the $\alpha$-cut $\Gamma^{\alpha}$ of $\Gamma$ as follows 
\begin{equation*}
\begin{split}
\Gamma^{\alpha} : \Omega \to cwk(X), 
\quad
\Gamma^{\alpha}(t) = [\Gamma(t)]^{\alpha} = \{x \in X : \Gamma(t)(x) \geq \alpha\}
\quad
(\alpha \in (0,1]).  
\end{split}
\end{equation*}
We say that the fuzzy map 
$\Gamma$ 
is scalarly measurable (scalarly integrable) if 
the map $\sigma(x^{*}, \Gamma^{\alpha}(\cdot))$ is measurable (Lebesgue integrable) 
for each $x^{*} \in X^{*}$ 
and $\alpha \in (0,1]$. 
The fuzzy map $\Gamma$ is said to be \textit{Pettis integrable} if 
\begin{itemize}
\item[(i)]
$\Gamma$ is scalarly integrable, 
\item[(ii)]
for each $A \in \Sigma$ there exists 
a generalized fuzzy number $M(A) \in \mathcal{F}_{cwk}(X)$ such that  
\begin{equation*}
\begin{split}
\sigma \left ( x^{*}, [M(A)]^{\alpha} \right ) 
= 
\int_{A} \sigma(x^{*}, \Gamma^{\alpha}(t)) d \mu
~\text{ for all }x^{*} \in X^{*}, \alpha \in (0,1].
\end{split}
\end{equation*}
\end{itemize}
We call $M(A)$ the \textit{fuzzy Pettis integral} of $\Gamma$ over $A$ and 
set $(FP)\int_{A} \Gamma(t) d \mu = M(A)$. 
\end{definition}

\section{The main results}

The main results are Theorem \ref{T_Pettis_WCG} and Theorem \ref{CT_Pettis_MultiF}. 
The convergence theorem for Pettis integrable multifunctions (Theorem \ref{CT_Pettis_MultiF}) 
together with  \cite[Lemma 5.4]{KAL_Pettis_LCS} yields 
a convergence theorem for Pettis integrable fuzzy maps, 
Theorem \ref{CT_Pettis_FuzzyMaps}. 
Let us start with a few auxiliary lemmas. 
The first lemma shows that each multimeasure with values in the hyperspace $cwk(X)$ 
is an $h$-multimeasure.

\begin{lemma}\label{L_Multi_HausdorffM}
If $M : \Sigma \to cwk(X)$ is a multimeasure, 
then $M$ is an $h$-multimeasure.
\end{lemma}
\begin{proof}
Let $(A_{n})$ be an arbitrary pairwise disjoint sequence of 
sets in $\Sigma$ with $A = \bigcup_{n} A_{n}$ 
and let $i \in \mathbb{N}$. 
Note that 
\begin{equation*}
\begin{split}
M_{i} : \Sigma \to cwk(X_{i}), 
\quad
M_{i}(E) = \varphi_{i}(M(E)),
\end{split}
\end{equation*}
where $cwk(X_{i})$ is the family of all nonempty convex and weakly compact subsets of $X_{i}$, 
and since 
\begin{equation*}
\begin{split}
H_{p_{i}}
\left ( 
\sum_{k=1}^{n} M(A_{k}), M(A)
\right )
=
h_{i}
\left ( 
\varphi_{i} \left ( \sum_{k=1}^{n} M(A_{k}) \right ), \varphi_{i}(M(A))
\right )
=
h_{i}
\left (
\sum_{k=1}^{n} M_{i}(A_{k}), M_{i}(A)
\right )
\end{split}
\end{equation*}
it is enough to prove that 
\begin{equation}\label{eqL_Multi_HausdorffM.1}
\begin{split}
\lim_{n \to \infty} 
h_{i}
\left (
\sum_{k=1}^{n} M_{i}(A_{k}), M_{i}(A)
\right )=0
\qquad
(\text{see Definition \ref{D_Multimeasures}}).
\end{split}
\end{equation}

Given any $n^{*}_{i} \in X^{*}_{i}$, 
$x^{*} = n^{*}_{i} \circ \varphi_{i} \in X^{*}$. 
Then, from the fact that $M$ is a multimeasure it follows that 
\begin{equation*}
\begin{split}
\lim_{n \to \infty} 
\sum_{k=1}^{n} \sigma(n_{i}^{*}, M_{i}(A_{k}))
=
\lim_{n \to \infty} 
\sum_{k=1}^{n} \sigma(x^{*}, M(A_{k}))
=\sum_{k=1}^{+\infty} \sigma(x^{*}, M(A_{k}))
=\sum_{k=1}^{+\infty} \sigma(n_{i}^{*}, M_{i}(A_{k})).
\end{split}
\end{equation*} 
This means that $M_{i}$ is a multimeasure 
and, therefore by \cite[Theorem 8.4.10, p.852]{HU} the map $M_{i}$ is 
an $h_{i}$-multimeasure. 
Hence, \eqref{eqL_Multi_HausdorffM.1} holds and the proof is complete.
\end{proof}


\begin{lemma}\label{L_Pettis_WCGII}
Let  $\Gamma : \Omega \to cwk(X)$ be a scalarly integrable multifunction. 
Assume that 
\begin{itemize}
\item
$\Gamma$ is determined by a WCG space $Y \subset X$, 
\item
$\varphi_{\Gamma} : X^{*} \to L_{1}(\mu)$ is weakly compact.
\end{itemize}
Then, for each $i \in \mathbb{N}$, we have
\begin{itemize}
\item[(i)]
$\Gamma_{i}(\cdot) = \varphi_{i}(\Gamma(\cdot))$ is determined by a WCG space $Y_{i} \subset X_{i}$, 
\item[(ii)]
$\varphi_{\Gamma_{i}} : X^{*}_{i} \to L_{1}(\mu)$ is weakly compact.
\end{itemize}
\end{lemma}
\begin{proof}
$(i)$ 
Since $Y$ is a WCG subspace of $X$ there is a weakly compact set $W \subset X$ such that 
$Y = \overline{\text{sp}(W)}$. 
According to \cite[Proposition 8, p.34]{ROB} 
(or, \cite[Corollary 2, p.65]{SCH}), 
the map  
$\varphi_{i}: (X,\sigma(X,X^{*})) \to (X_{i},\sigma(X_{i},X_{i}^{*}))$ is continuous. 
Then,  
$\varphi_{i}(W)$ is a weakly compact subset of $X_{i}$ whenever $i \in \mathbb{N}$. 
Thus, 
\begin{equation*} 
\begin{split}
Y_{i} 
= \overline{\text{sp}(\varphi_{i}(W))}^{||\cdot||_{i}} 
= \overline{\varphi_{i}(\text{sp}(W))}^{||\cdot||_{i}}
\qquad
\left (
\overline{\varphi_{i}(\text{sp}(W))}^{||\cdot||_{i}}
\supset
\varphi_{i} \left ( \overline{\text{sp}(W)} \right )
\right ) 
\end{split}
\end{equation*}
is a WCG subspace of $X_{i}$ 
and hence 
\begin{equation*}
\begin{split}
n^{*}_{i} \in Y_{i}^{\perp}
\Rightarrow& 
[n^{*}_{i}(\varphi_{i}(x)) =0\text{ for all }x \in \overline{\text{sp}(W)}=Y]
\Rightarrow
n^{*}_{i} \circ \varphi_{i} \in Y^{\perp}
\\
\Rightarrow&
\sigma(n^{*}_{i}, \Gamma_{i}(t)) = \sigma(n^{*}_{i} \circ \varphi_{i}, \Gamma(t)) = 0
~\text{ at almost all }t \in \Omega.
\end{split}
\end{equation*}
Consequently, $\Gamma_{i}$ is determined by the WCG space $Y_{i}$.

$(ii)$ 
Since
\begin{equation*}
\begin{split}
x^{*} \in U_{p_{i}}^{o} 
\Leftrightarrow 
(\exists n^{*}_{i} \in B_{X_{i}^{*}})
[x^{*} = n^{*}_{i} \circ \varphi_{i}]
\qquad
(i \in \mathbb{N}), 
\end{split}
\end{equation*} 
it follows that 
\begin{equation*}
\begin{split}
\{\sigma(n_{i}^{*}, \Gamma_{i}(\cdot)): n_{i}^{*} \in B_{X_{i}^{*}}  \}
=
\{\sigma(x^{*}, \Gamma(\cdot)): x^{*} \in U_{p_{i}}^{o}\},
\qquad
(i \in \mathbb{N}). 
\end{split}
\end{equation*}
Hence, 
$\{\sigma(n_{i}^{*}, \Gamma_{i}(\cdot)): n_{i}^{*} \in B_{X_{i}^{*}}  \}$ 
is a relatively weakly compact subset of $L_{1}(\mu)$, 
and this completes the proof.
\end{proof}

Given $i \leq j$ define a continuous linear map 
$g_{ij} : X_{j} \to X_{i}$ by $g_{ij}(\varphi_{j}(x)) = \varphi_{i}(x)$ for all $x \in X$, 
and let $\overline{g}_{ij} : \overline{X}_{j} \to \overline{X}_{i}$ be the continuous extension of $g_{ij}$ to 
$\overline{X}_{j}$. 
Thanks to \cite[5.4, p.53]{SCH}, 
the projective limit of the normed spaces $(X_{i})$ with respect to $(g_{ij})_{i \leq j}$ 
is equal to  
the projective limit of the Banach spaces $(\overline{X}_{i})$ with respect to $(\overline{g}_{ij})_{i \leq j}$, 
and the map
\begin{equation}\label{eqProjective_Limit}
\begin{split}
\phi: X \to Y = \underset{\longleftarrow}{\lim} g_{ij} X_{j} = \underset{\longleftarrow}{\lim} \overline{g}_{ij} \overline{X}_{j}, 
\quad
\phi(x) = (\varphi_{i}(x)),
\end{split}
\end{equation}
is an isomorphism.


The following lemma is the key to prove 
the main results, Theorem \ref{T_Pettis_WCG} and Theorem \ref{CT_Pettis_MultiF}.

\begin{lemma}\label{L_Key_WCG}
Let  $\Gamma : \Omega \to cwk(X)$ be a scalarly integrable multifunction such that 
the multifunction $\Gamma_{i} = \varphi_{i}(\Gamma(\cdot))$ is Pettis integrable in 
$cwk(\overline{X}_{i})$
with $M_{i}(A) = (P)\int_{A} \Gamma_{i}(t) d \mu$ for all $i \in \mathbb{N}$ and $A \in \Sigma$. 

Then, $\Gamma$ is Pettis integrable in $cwk(X)$ with 
$M(A) = (P)\int_{A} \Gamma(t) d \mu$ for all $A \in \Sigma$, and 
\begin{equation}\label{eqL_Pettis_WCGIII.0}
\begin{split}
\varphi_{i}(M(A)) = M_{i}(A)
~\text{ for every }A \in \Sigma, i \in \mathbb{N}. 
\end{split}
\end{equation}
\end{lemma}
\begin{proof}
Let $m^{*}_{i}: \overline{X}_{i} \to \mathbb{R}$ be a continuous linear functional 
and let $i \leq j$. 
Then, for each $A \in \Sigma$, 
\begin{equation*}
\begin{split}
\sigma(m^{*}_{i}, M_{i}(A))
=& 
\int_{A} \sigma(m^{*}_{i}, \Gamma_{i}(t)) d \mu
=
\int_{A} \sigma(m^{*}_{i}, \overline{g}_{ij} (\Gamma_{j}(t))) d \mu
=
\int_{A} \sigma(m^{*}_{i} \circ \overline{g}_{ij}, \Gamma_{j}(t)) d \mu
\\
=&
\sigma(m^{*}_{i} \circ \overline{g}_{ij}, M_{j}(A))
=
\sigma(m^{*}_{i}, \overline{g}_{ij} (M_{j}(A))).
\end{split}
\end{equation*}
Thus, 
\begin{equation*}
\begin{split}
\sigma(m^{*}_{i}, M_{i}(A))
=
\sigma(m^{*}_{i}, \overline{g}_{ij} (M_{j}(A)))
~\text{ for every }m^{*}_{i} \in \overline{X}_{i}^{*}, \quad (A \in \Sigma), 
\end{split}
\end{equation*}
where $\overline{X}_{i}^{*}$ is the topological dual of the Banach space $\overline{X}_{i}$. 
The last result together with the Hahn-Banach separation
theorem (\cite[Theorem V.2.10]{DUN}) yields that
\begin{equation}\label{eqL_Pettis_WCGIII.1}
\begin{split}
M_{i}(A)
=
\overline{g}_{ij} (M_{j}(A))
\qquad
(A \in \Sigma, i \leq j). 
\end{split}
\end{equation}

$\bullet$ 
We will prove that  
\begin{equation*}
\begin{split}
M_{i}(A) \subset X_{i}
\qquad
(i \in \mathbb{N}, A \in \Sigma). 
\end{split}
\end{equation*}
To do this, let $w^{0}_{i}$ be a member  of $M_{i}(A)$. 
For $s=i, i+1, i+2, \dotsc$,  inductively  there exists $w^{0}_{s+1} \in M_{s+1}(A)$ such that 
$\overline{g}_{s,s+1}(w^{0}_{s+1}) = w^{0}_{s}$ (see \eqref{eqL_Pettis_WCGIII.1}), 
while if $1<s\leq i$ define 
\begin{equation*}
\begin{split}
w^{0}_{i-1} = \overline{g}_{i-1,i}(w^{0}_{i}),
~
w^{0}_{i-2} = \overline{g}_{i-2,i-1}(w^{0}_{i-1}),
\dotsc, 
~
w^{0}_{1} = \overline{g}_{1,2}(w^{0}_{2})
\end{split}
\end{equation*}
Thus,  
the vector 
$w^{0} = (w^{0}_{s}) \in \prod_{s=1}^{+\infty} M_{s}(A) 
\subset \prod_{s=1}^{+\infty} \overline{X}_{i}$ is such that 
\begin{equation*}
\begin{split}
\overline{g}_{s,s+1}(w^{0}_{s+1}) = w^{0}_{s}
~\text{ for every }s \in \mathbb{N}.
\end{split}
\end{equation*}
If we write $\overline{g}_{s,s+1} = h_{s+1}$ for all  $s \in \mathbb{N}$, then  
\begin{equation*}
\begin{split}
\overline{g}_{s, s + k} = h_{s+1} \circ  h_{s+2} \circ \dotsc \circ  h_{s+k}
\qquad
(s, k \in \mathbb{N}).
\end{split}
\end{equation*}
Hence,  
\begin{equation*}
\begin{split}
w^{0} = (w^{0}_{s}) \in \underset{\longleftarrow}{\lim} \overline{g}_{ij} \overline{X}_{j}
\end{split}
\end{equation*}
and, therefore, there exists $x_{0} \in X$ such that $\phi(x_{0}) = w^{0}$ (see \eqref{eqProjective_Limit}). 
Since $\varphi_{i}(x_{0}) = w^{0}_{i} \in X_{i}$, 
we infer that $M_{i}(A) \subset X_{i}$ ($i \in \mathbb{N}$, $A \in \Sigma$). 
Consequently, 
\begin{equation}\label{eqL_Pettis_WCGIII.11}
\begin{split}
M_{i}(A)
=
g_{ij} (M_{j}(A))
\quad
(A \in \Sigma, i \leq j)
\quad
(\text{see \eqref{eqL_Pettis_WCGIII.1}}).
\end{split}
\end{equation}

$\bullet$ 
Since
\begin{equation*}
\begin{split}
\left \{ (w_{i}) \in \prod_{i} M_{i}(A) : i \leq j \Rightarrow g_{ij}(w_{j}) = w_{i} \right \}
\subset \underset{\longleftarrow}{\lim} g_{ij} X_{j}
\end{split}
\end{equation*}
there is $M(A) \subset X$ such that 
\begin{equation}\label{eqL_Pettis_WCGIII.02}
\begin{split}
\phi(M(A)) = \left \{ (w_{i}) \in \prod_{i} M_{i}(A) : i \leq j \Rightarrow g_{ij}(w_{j}) = w_{i} \right \}. 
\end{split}
\end{equation}
We are going to prove that 
\begin{equation}\label{eqL_Pettis_WCGIII.2}
\begin{split}
\varphi_{i}(M(A)) = M_{i}(A)
\qquad(i \in \mathbb{N}).
\end{split}
\end{equation}
Since $\varphi_{i}(M(A)) \subset M_{i}(A)$, it remains to prove the inverse inclusion. 
To do this, let $u^{0}_{i}$ be a member  of $M_{i}(A)$. 
For $s=i, i+1, i+2, \dotsc$,  inductively  there exists $u^{0}_{s+1} \in M_{s+1}(A)$ such that 
$g_{s,s+1}(u^{0}_{s+1}) = u^{0}_{s}$ (see \eqref{eqL_Pettis_WCGIII.11}), 
while if $1<s \leq i$ define 
\begin{equation*}
\begin{split}
u^{0}_{i-1} = g_{i-1,i}(u^{0}_{i}), 
~u^{0}_{i-2} = g_{i-2,i-1}(u^{0}_{i-1}),
\dotsc,
~u^{0}_{1} = g_{1,2}(u^{0}_{2}).
\end{split}
\end{equation*}
Thus, the vector 
$u^{0} = (u^{0}_{s}) \in \prod_{s=1}^{+\infty} M_{s}(A) 
\subset \prod_{s=1}^{+\infty} X_{i}$ is such that 
\begin{equation*}
\begin{split}
g_{s,s+1}(u^{0}_{s+1}) = u^{0}_{s}
~\text{ for every }s \in \mathbb{N}.
\end{split}
\end{equation*}
Just as above we show that 
\begin{equation*}
\begin{split}
g_{s,s+k}(u^{0}_{s+k}) = u^{0}_{s}
~\text{ for every }s,k \in \mathbb{N}.
\end{split}
\end{equation*}
The last result together with \eqref{eqL_Pettis_WCGIII.02} yields that 
$u^{0} = (u^{0}_{s}) \in \phi(M(A))$. 
Hence,  
$u^{0}_{i} \in \varphi_{i}(M(A))$. 
Consequently, $M_{i}(A) \subset \varphi_{i}(M(A))$, 
and since $\varphi_{i}(M(A)) \subset M_{i}(A)$, 
it follows that \eqref{eqL_Pettis_WCGIII.2} holds.

$\bullet$ 
Given $x^{*}_{i} \in X^{*}_{p_{i}}$   
there exists $n^{*}_{i} \in X^{*}_{i}$ such that $x^{*}_{i} = n^{*}_{i} \circ \varphi_{i}$. 
Hence, for any $A \in \Sigma$ we obtain by \eqref{eqL_Pettis_WCGIII.2} that
\begin{equation*}
\begin{split}
\sigma(x^{*}_{i}, M(A)) 
=& 
\sigma(n^{*}_{i} \circ \varphi_{i}, M(A))
=
\sigma(n^{*}_{i}, M_{i}(A)) 
= \int_{A} \sigma(n^{*}_{i}, \Gamma_{i}(t)) d \mu 
\\
=&
\int_{A} \sigma(n^{*}_{i} \circ \varphi_{i}, \Gamma(t)) d \mu 
=
\int_{A} \sigma(x^{*}_{i}, \Gamma(t)) d \mu,
\end{split}
\end{equation*}
and since 
$X^{*} = \bigcup_{i=1}^{+\infty} X^{*}_{p_{i}}$, 
it follows that $\Gamma$ is Pettis integrable and \eqref{eqL_Pettis_WCGIII.0} holds,
and this completes the proof.
\end{proof}

We are now ready to present a full characterization of Pettis integrability of a 
scalarly integrable multifunction in terms of a WCG subspace of $X$.

\begin{theorem}\label{T_Pettis_WCG}
Let  $\Gamma : \Omega \to cwk(X)$ be a scalarly integrable multifunction. 
Then the following statements are equivalent:
\begin{itemize}
\item[(i)]
$\Gamma$ is Pettis integrable in $cwk(X)$, 
\item[(ii)]
$\Gamma$ is determined by a WCG space $Y \subset X$ 
and $\varphi_{\Gamma} : X^{*} \to L_{1}(\mu)$ is weakly compact. 
\end{itemize}
\end{theorem}
\begin{proof}
$(i) \Rightarrow (ii)$ 
Assume that $\Gamma$ is Pettis integrable in $cwk(X)$. 
According to Theorem \ref{T_CWK} 
the set $M(\Sigma) = \bigcup_{E \in \Sigma} M(E)$ is a relatively weakly compact subset of $X$. 
Then, $Y = \overline{\text{sp}(M(\Sigma))}$ is a WCG subspace of $X$. 
Hence, for each $E \in \Sigma$ we have 
\begin{equation*}
\begin{split}
x^{*} \in Y^{\perp}
&\Rightarrow 
[(\forall y \in Y)[x^{*}(y)=0]\text{ and }Y \supset M(E)] 
\Rightarrow
\sigma(x^{*}, M(E)) = 0
\\
&\Rightarrow
0 = \sigma(x^{*}, M(E)) = \int_{E} \sigma(x^{*}, \Gamma(t)) d \mu.
\end{split}
\end{equation*}
Thus, 
\begin{equation*}
\begin{split}
x^{*} \in Y^{\perp} \Rightarrow 
(\forall E \in \Sigma)
\left [
0 = \int_{E} \sigma(x^{*}, \Gamma(t)) d \mu
\right ]
\end{split}
\end{equation*}
and, consequently, 
\begin{equation*}
\begin{split}
x^{*} \in Y^{\perp} \Rightarrow 
\sigma(x^{*}, \Gamma(t)) =0
~\text{ at almost all }t \in \Omega.
\end{split}
\end{equation*}
The last result together with Definition \ref{D_DeterminedMultiF_By_WCG} yields that 
$\Gamma$ is determined by the WCG space $Y$.

It remains to prove that $\varphi_{\Gamma} : X^{*} \to L_{1}(\mu)$ is weakly compact. 
To do this, let $(A_{n})$ be an arbitrary pairwise disjoint sequence of sets in $\Sigma$ with $A = \bigcup_{n} A_{n}$ 
and let $i \in \mathbb{N}$. 
Since the multimeasure $M: \Sigma \to cwk(X)$ is a strong multimeasure (or, equivalently multimeasure), 
it follows by Lemma \ref{L_Multi_HausdorffM} that $M$ is an $h$-multimeasure. 
Then,
\begin{equation*}
\begin{split}
\lim_{n \to \infty} 
H_{p_{i}} \left ( \sum_{k=1}^{n} M(A_{k}), M(A) \right ) = 0
\end{split}
\end{equation*}
and since we have by \cite[Theorem II-18, p.49]{CAST} that
\begin{equation*}
\begin{split}
H_{p_{i}} \left ( \sum_{k=1}^{n} M(A_{k}), M(A) \right )
=&
\sup_{x^{*} \in U_{p_{i}}^{o}}
\left |
\sigma \left (x^{*}, \sum_{k=1}^{n} M(A_{k}) \right ) 
-
\sigma \left (x^{*}, M(A) \right ) 
\right |
\\
=&
\sup_{x^{*} \in U_{p_{i}}^{o}}
\left |
 \sum_{k=1}^{n} \sigma(x^{*}, M(A_{k})) 
-
\sigma \left (x^{*}, M(A) \right ) 
\right |
\\
=&
\sup_{x^{*} \in U_{p_{i}}^{o}}
\left |
\sum_{k > n} \sigma(x^{*},  M(A_{k}) )  
\right |,
\end{split}
\end{equation*}
it follows that
\begin{equation*}
\begin{split} 
\lim_{n \to \infty} \sum_{k > n} \sigma(x^{*}, M(A_{k})) = 0
~\text{ uniformly for }x^{*} \in U_{p_{i}}^{o}.
\end{split}
\end{equation*}
Since $(A_{n})$ was arbitrary,   
the last result shows that 
$\{\sigma(x^{*}, M(\cdot)): x^{*} \in U_{p_{i}}^{o}\}$ 
is an uniformly countable additive family of real-valued measures 
(cf. \cite[Definition 14, p.7]{DIES}), 
and since $M(\Sigma)$ is bounded in $X$ (see \cite[Lemma 4.5]{KAL_Pettis_LCS}),  
it follows by \cite[Corollary 5, p.13]{DIES} that 
\begin{equation*}
\begin{split}
0=
\lim_{\mu(E) \to 0} \sigma(x^{*}, M(E))
=
\lim_{\mu(E) \to 0} \left |\int_{E} \sigma(x^{*}, \Gamma(t)) d \mu \right |
~\text{ uniformly for }x^{*} \in U_{p_{i}}^{o}.
\end{split}
\end{equation*}
Further, we obtain by \cite[Corollary IV.8.11, p.294]{DUN} that 
$\{\sigma(x^{*}, \Gamma(\cdot)): x^{*} \in U_{p_{i}}^{o}\}$  
is weakly sequentially compact subset of $L_{1}(\mu)$ 
and hence by \cite[Theorem V.6.1, p.430]{DUN},  
it is relatively weakly compact subset of $L_{1}(\mu)$. 
This means that $\varphi_{\Gamma}$ is weakly compact.

$(ii) \Rightarrow (i)$ 
Assume that $(ii)$ holds. 
Then, we obtain by Lemma \ref{L_Pettis_WCGII} that 
$\Gamma_{i}(\cdot) = \varphi_{i}(\Gamma(\cdot))$ is determined by a WCG space 
$Y_{i} \subset X_{i}$ and 
the map $\varphi_{\Gamma_{i}} : X^{*}_{i} \to L_{1}(\mu)$ is weakly compact 
($i \in \mathbb{N}$).
Hence, by \cite[Theorem 2.5]{MUS_Pettis_MultiF} 
the multifunction $\Gamma_{i}$ is Pettis integrable in $cwk(\overline{X}_{i})$ with 
$M_{i}(A) = (P)\int_{A} \Gamma(t) d \mu$ for all $A \in \Sigma$ ($i \in \mathbb{N}$). 
Further, by Lemma \ref{L_Key_WCG} the multifunction 
$\Gamma$ is Pettis integrable in $cwk(X)$ and the proof is complete.
\end{proof}


The following is a convergence theorem for Pettis integrable multifunctions.

\begin{theorem}\label{CT_Pettis_MultiF}
Let $\Gamma : \Omega \to cwk(X)$ be a  scalarly integrable multifunction
and let $\{\Gamma_{n} : \Omega \to cwk(X) : n \in \mathbb{N}\}$ be a sequence of multifunctions 
Pettis integrable in $cwk(X)$ and satisfying two conditions:
\begin{itemize}
\item[(a)]
the set 
$\left \{\sigma(x^{*}, \Gamma_{n}(\cdot)) : x^{*} \in U_{p_{i}}^{o}, n \in \mathbb{N} \right \}$ 
is uniformly integrable whenever $i \in \mathbb{N}$, 
\item[(b)]
for each $x^{*} \in X^{*}$, 
\begin{equation*}
\begin{split}
\lim_{n \to \infty} \sigma(x^{*}, \Gamma_{n}) = \sigma(x^{*}, \Gamma) 
~\text{in }\mu\text{-measure}.
\end{split}
\end{equation*}
\end{itemize} 
Then $\Gamma$ is Pettis integrable in $cwk(X)$ and, 
\begin{equation}\label{eqC_Pettis_WCG.0}
\begin{split}
\lim_{n \to \infty} 
\sigma \left (x^{*}, (P)\int_{A} \Gamma_{n}(t) d \mu \right ) 
=
\sigma \left (x^{*}, (P)\int_{A} \Gamma(t) d \mu \right ) 
~\text{ for every }x^{*} \in X^{*}, A \in \Sigma.
\end{split}
\end{equation}
\end{theorem}
\begin{proof}
Fix an arbitrary $i \in \mathbb{N}$ and let $\overline{X}^{*}_{i}$ be the topological dual of the Banach space 
$\overline{X}_{i}$. 
We write $\Gamma_{i}(\cdot) = \varphi_{i}(\Gamma(\cdot))$ and $\Gamma^{i}_{n}(\cdot) = \varphi_{i}(\Gamma_{n}(\cdot))$ for every 
$n \in \mathbb{N}$. 
Since $\sigma(n_{i}^{*} \circ \varphi_{i}, \Gamma_{n}(\cdot)) = \sigma(n_{i}^{*}, \Gamma^{i}_{n}(\cdot))$ for all 
$n_{i}^{*} \in X^{*}_{i}$  and $n \in \mathbb{N}$, 
it follows that 
\begin{equation}\label{eqC_Pettis_WCG.1} 
\begin{split}
\int_{A} \sigma(n_{i}^{*}, \Gamma^{i}_{n}(t)) d \mu 
=&
\int_{A} \sigma(n_{i}^{*} \circ \varphi_{i}, \Gamma_{n}(t)) d \mu 
=
\sigma \left (n_{i}^{*} \circ \varphi_{i}, (P)\int_{A} \Gamma_{n}(t) d \mu \right )
\\
=&
\sigma \left ( n_{i}^{*}, \varphi_{i} \left ( (P)\int_{A} \Gamma_{n}(t) d \mu \right ) \right )
\quad
(A \in \Sigma).
\end{split}
\end{equation}
This means that each multifunction $\Gamma_{n}^{i}$ is Pettis integrable in $cwk(X_{i})$. 
Since $n_{i}^{*} \circ \varphi_{i} \in U_{p_{i}}^{o}$ whenever  $n_{i}^{*} \in B_{X^{*}_{i}}$, we have also 
\begin{itemize}
\item
the set 
$\{\sigma(n_{i}^{*}, \Gamma^{i}_{n}(\cdot)) : n_{i}^{*} \in B_{X^{*}_{i}}, n \in \mathbb{N}\}$ 
is uniformly integrable, 
\item
for each $n_{i}^{*} \in X_{i}^{*}$, the sequence $(\sigma(n_{i}^{*}, \Gamma^{i}_{n}))$ 
converges in measure to $\sigma(n_{i}^{*}, \Gamma_{i})$.
\end{itemize}  
Therefore, by \cite[Theorem 5.1]{MUS_Pettis_MultiF}, 
the multifunction $\Gamma_{i}$ 
is Pettis integrable in $cwk \left (\overline{X}_{i} \right )$ and
\begin{equation}\label{eqC_Pettis_WCG.2}  
\begin{split}
\lim_{n \to \infty} 
\sigma \left (m^{*}_{i}, (P)\int_{A} \Gamma_{n}^{i}(t) d \mu \right ) 
=
\sigma \left (m^{*}_{i}, (P)\int_{A} \Gamma_{i}(t) d \mu \right ) 
~\text{ for every }m^{*}_{i} \in \overline{X}^{*}_{i}, A \in \Sigma.
\end{split}
\end{equation}
Since $i$ was arbitrary, the last result holds for every $i \in \mathbb{N}$. 
Hence, by Lemma \ref{L_Key_WCG}, the multifunction $\Gamma$ is Pettis integrable in $cwk(X)$ and
\begin{equation*}
\begin{split}
(P)\int_{A} \Gamma_{i}(t) d \mu
=
\varphi_{i} \left ( (P) \int_{A} \Gamma(t) d \mu \right )
\qquad
(i \in \mathbb{N}, A \in \Sigma) 
\end{split}
\end{equation*}
and by \eqref{eqC_Pettis_WCG.1}, 
\begin{equation*}
\begin{split}
(P)\int_{A} \Gamma_{n}^{i}(t) d \mu = \varphi_{i} \left ( (P) \int_{A} \Gamma_{n}(t) d \mu \right )
\qquad
(A \in \Sigma, n, i \in \mathbb{N}).
\end{split}
\end{equation*}
The last result together with \eqref{eqC_Pettis_WCG.2} yields
\begin{equation*}
\begin{split}
\lim_{n \to \infty} 
\sigma \left (n^{*}_{i}, (P)\int_{A} \Gamma_{n}^{i}(t) d \mu \right ) 
=
\sigma \left (n^{*}_{i}, (P)\int_{A} \Gamma_{i}(t) d \mu \right ) 
\qquad
(n^{*}_{i} \in X^{*}_{i}, A \in \Sigma).
\end{split}
\end{equation*}
Hence, 
for any $x^{*} \in X^{*}_{p_{i}}$, 
\begin{equation*}
\begin{split}
\lim_{n \to \infty} 
\sigma \left (x^{*}, \int_{A} \Gamma_{n}(t) d \mu \right )
=&
\lim_{n \to \infty} 
\sigma \left (n_{i}^{*} \circ \varphi_{i}, (P)\int_{A} \Gamma_{n}(t) d \mu \right )
=
\lim_{n \to \infty} 
\sigma \left (n_{i}^{*}, (P) \int_{A} \Gamma_{n}^{i}(t) d \mu \right )
\\
=&\sigma \left (n_{i}^{*}, (P) \int_{A} \Gamma_{i}(t) d \mu \right )
=\sigma \left (n_{i}^{*} \circ \varphi_{i}, (P) \int_{A} \Gamma(t) d \mu \right )
\\
=&
\sigma \left (x^{*}, (P)\int_{A} \Gamma(t) d \mu \right ),
\end{split}
\end{equation*} 
where  $n_{i}^{*} \in X^{*}_{i}$ is such that $x^{*} = n_{i}^{*} \circ \varphi_{i}$. 
Since $X^{*} = \bigcup_{i=1}^{+\infty} X^{*}_{p_{i}}$, the last result yields that 
\eqref{eqC_Pettis_WCG.0} holds, and this completes the proof.
\end{proof}

The next is a convergence theorem for Pettis integrable fuzzy maps.

\begin{theorem}\label{CT_Pettis_FuzzyMaps}
Let $\Gamma : \Omega \to \mathcal{F}_{cwk}(X)$ be a scalarly integrable fuzzy map
and let $\{\Gamma_{n} : \Omega \to \mathcal{F}_{cwk}(X) : n \in \mathbb{N}\}$ be a sequence of 
Pettis integrable fuzzy maps satisfying two conditions:
\begin{itemize}
\item[(a)]
the set 
$\{\sigma(x^{*}, \Gamma_{n}^{\alpha}(\cdot)) : x^{*} \in U_{p_{i}}^{o}, n \in \mathbb{N}\}$ 
is uniformly integrable whenever $i \in \mathbb{N}$ and $\alpha \in (0,1]$, 
\item[(b)]
for each $x^{*} \in X^{*}$ and $\alpha \in (0,1]$, 
\begin{equation}\label{eqCT_Pettis_FuzzyMaps.0}
\begin{split}
\lim_{n \to \infty} \sigma(x^{*}, \Gamma^{\alpha}_{n}) = \sigma(x^{*}, \Gamma^{\alpha}) 
~\text{in }\mu\text{-measure}.
\end{split}
\end{equation}
\end{itemize} 
Then $\Gamma$ is Pettis integrable  and, 
\begin{equation*}
\begin{split}
\lim_{n \to \infty} 
\sigma \left (x^{*}, (P)\int_{A} \Gamma^{\alpha}_{n}(t) d \mu \right ) 
=
\sigma \left (x^{*}, (P)\int_{A} \Gamma^{\alpha}(t) d \mu \right ) 
~\text{ for every }x^{*} \in X^{*}, A \in \Sigma, \alpha \in (0,1].
\end{split}
\end{equation*}
\end{theorem}
\begin{proof}
For each $\alpha \in (0,1]$, 
by Theorem \ref{CT_Pettis_MultiF} the map
\begin{equation*}
\begin{split}
\Gamma^{\alpha} : \Omega \to cwk{X}, 
\quad
\Gamma^{\alpha}(t) = [\Gamma(t)]^{\alpha},
\quad
[\Gamma(t)]^{\alpha} = \{x \in X : \Gamma(t)(x) \geq \alpha\},
\end{split}
\end{equation*}
is Pettis integrable and \eqref{eqCT_Pettis_FuzzyMaps.0} holds. 
Hence, we obtain by \cite[Lemma 5.4]{KAL_Pettis_LCS} that the fuzzy map $\Gamma$ is Pettis integrable 
and the proof is complete.
\end{proof}


\begin{thebibliography}{plain}





\bibitem{BAN} J.~Ban, 
\textit{Radon-Nikod\'{y}m theorem and conditional expectation 
of fuzzy valued measures and variables}, 
Fuzzy Sets and Systems \textbf{34} (1990), 383-392







\bibitem{CASC_SelectorsA} B.~Cascales, V.~Kadets, and J.~Rodr\'{i}guez,  
\textit{Measurable selectors and set-valued Pettis integral in non-separable Banach spaces}, 
J. Funct. Anal. \textbf{256/3} (2009), 673-699.






\bibitem{CASC_SelectorsB} B.~Cascales, V.~Kadets, and J.~Rodr\'{i}guez, 
\textit{Measurability and selections of multifunctions in Banach spaces}, 
J. Convex Analysis \textbf{17} (2010), 229-240. 







\bibitem{CAST} C. ~Castaing, ~M. ~Valadier, 
\textit{Convex Analysis and Measurable Multifunctions} ,
Lect. Notes Math., \textbf{580}, Springer-Verlag, Berlin (1977).





\bibitem{DIES} J.~Diestel and J.~J.~Uhl, 
\textit{Vector Measures}, 
Math. Surveys, vol.15, Amer. Math. Soc., Providence, RI, 1977.




\bibitem{DUN} N. Dunford and J. T. Schwartz, 
\textit{Linear operators. Part I: General theory}, 
Interscience, New York, (1958).







\bibitem{EL} K.~El Amri and C.~Hess, 
\textit{On the Pettis integral of closed valued multifunctions}, 
Set-Valued Anal. \textbf{8(4)} (2000), 329-360.


\bibitem{HU} S.~Hu, N.~S.~Papageorgiou, 
\textit{Handbook of Multivalued Analysis I}, 
Kluwer, Dordrecht (1997).



\bibitem{JAMES} R.~C.~James, 
\textit{Weakly compact sets}, 
Trans. Amer. Math. Soc. \textbf{113} (1964), 129-140.






\bibitem{KAL_REND} S.~B.~Kaliaj, 
\textit{Indefinite Pettis integral of multifunctions  in locally convex spaces}, 
Rend. Istit. Mat. Univ. Trieste \textbf{53} (2020). 






\bibitem{KAL_GFNM} S.~B.~Kaliaj, 
\textit{Radon-Nikod\'{y}m theorems for generalized fuzzy number measures}, 
Fuzzy Sets and Systems \textbf{512} (2025).




\bibitem{KAL_FNM} S.~B.~Kaliaj, 
\textit{Radon-Nikod\'{y}m theorems for fuzzy number multimeasures}, 
Fuzzy Sets and Systems \textbf{524} (2026).






\bibitem{KAL_FNM_Frechet} S.~B.~Kaliaj, 
\textit{A Radon-Nikod\'{y}m theorem for multimeasures in Fr\'{e}chet spaces with application}, 
J. Math. Anal. Appl. \textbf{557/2} (2026).



\bibitem{KAL_Pettis_LCS} S.~B.~Kaliaj, 
\textit{Characterizations of Pettis Integrable Multifunctions in locally convex spaces} .
arXiv:2611.16969 [math.GM] (2026), 
https://doi.org/10.48550/arXiv.2611.16969.








\bibitem{KAN} D.~A.~Kandilakis, 
\textit{On the extension of multimeasures and integration with respect to a multimeasure}, 
PAMS \textbf{116} (1992), 85-92.






\bibitem{KELL} J.~L.~Kelly,
\textit{General Topology}, 
Springer New York Van Nostrand, (1955).






\bibitem{KOTHE} G.~K\"{o}the, 
\textit{Topological Vector Spaces. I}, 
Grundlehren Math.~Wiss., vol.\textbf{159}, Springer, New York, (1969). 







\bibitem{MUS_Pettis_Integral_Functions} K.~Musial,
\textit{Topics in the theory of Pettis integration},
Rend. Istit. Mat. Univ. Trieste \textbf{23} (1991), 177-262. 







\bibitem{MUS_Pettis_MultiF} K.~Musial,
\textit{Pettis integrability of multifunctions with values in arbitrary Banach spaces}, 
J. Convex Analysis \textbf{18} (2011), 769-810.






\bibitem{MUS_Approximation} K.~Musial,
\textit{Approximation of Pettis integrable multifunctions with values in arbitrary Banach spaces}, 
J. Convex Analysis \textbf{20(3)} (2013), 833-870.


\bibitem{MUS_Conjugate_Banach_Spaces} K.~Musial,
\textit{Multimeasures with values in conjugate Banach spaces and the weak
Radon-Nikodým property}, 
J. Convex Analysis \textbf{28/3} (2021), 879-902.



\bibitem{MUS_RNT_Rich_Multimeasures_I} K.~Musial, 
\textit{Lifting Approach to Integral Representations of 
Rich Multimeasures with Values in Banach Spaces $I$},
J.~Convex Anal., \textbf{(31/1)} (2024), 179-193. 





\bibitem{MUS_RNT_Rich_Multimeasures_II} K.~Musial, 
\textit{Lifting Approach to Integral Representations of Rich Multimeasures with Values in Banach Spaces $II$},
J.~Convex Anal., \textbf{(31/3)} (2024), 867-888. 






\bibitem{ROB} A.~P.~Robertson and W.~Robertson, 
\textit{Topological Vector Spaces}, 
Cambridge University  Press (1964). 






\bibitem{SCH} H.~H.~Schaefer, W.~P.~Wolff, 
\textit{Topological Vector Spaces},
Originally published by Springer-Verlag New York, (1999).






\bibitem{STOJ_RNT} M.~Stojakovi\'{c}, 
\textit{Fuzzy valued measure}, 
Fuzzy Sets and Systems \textbf{65/1} (1994), 95-104.







\bibitem{TAL}  M.~Talagrand, 
\textit{Pettis Integral and Measure Theory}, 
Mem. Amer. Math. Soc. \textbf{51} (1984).





\bibitem{WU_RNT} Jianrong Wu, Xiaoping Xue, Congxin Wu, 
\textit{Radon-Nikod\'{y}m theorem and Vitali-Hahn-Saks theorem on fuzzy number measures in Banach spaces}, 
Fuzzy Sets and Systems \textbf{117/3} (2001), 339-346. 





\bibitem{ZAD} L.~A.~Zadeh, 
\textit{Fuzzy sets}, 
Inform. and Control \textbf{8} (1965), 338-353.





\bibitem{ZIAT} H.~Ziat, 
\textit{On a characterization of Pettis integrable multifunctions}, 
Bull. Polish Acad. Sci. Math. 48 (2000), 227-230.



\bibitem{ZHOU} Cai-Li Zhou and Fu-Gui Shi,
\textit{Lebesgue Decomposition Theorem and Weak Radon-Nikodým
Theorem for Generalized Fuzzy Number Measures}, 
Journal of Function Spaces Volume 2015, Article ID 576134, 8 pages,
http://dx.doi.org/10.1155/2015/576134.






\end{thebibliography}
\end{document}